\documentclass[12pt]{article}
\usepackage{amsfonts}
\usepackage{amsmath}
\usepackage[english]{babel}
\newtheorem{prop}{Proposition}[section]

\newtheorem{thm}[prop]{Theorem}
\newtheorem{lem}[prop]{Lemma}

\newtheorem{rem}[prop]{Remark}

\def\1{1\!{\rm l}}

\begin{document}
\title{The spectrum minimum for random Schr\"odinger operators with indefinite sign potentials} \maketitle
\begin{center}\bf{{Hatem Najar}} \footnote{D\'epartement de
Math\'ematiques Physiques I.P.E.I. Monastir, 5000 Monastir Tunisie
. \\Researches partially supported by CMCU N 04/S1404 and Research
Unity 01/UR/ 15-01 projects. }
\end{center}
\begin{abstract}

This paper sets out to study the spectral minimum for operator
belonging to the family of random Schr\"odinger operators of the
form $H_{\lambda,\omega}=-\Delta+W_{\text{per}}+\lambda
V_{\omega}$, where we suppose that $V_{\omega}$ is of Anderson type
and the single site is assumed to be with an indefinite sign.
Under some assumptions we prove that there exists $\lambda_0>0$
such that for any $\lambda \in [0,\lambda_0]$, the minimum of the
spectrum of $H_{\lambda,\omega}$ is obtained by a given
realization of the random variables. \vskip.5cm\noindent

\end{abstract}
{\small\sf 2000 Mathematics Subject Classification:81Q10, 35P05,
37A30,47F05.\\
Keywords and phrases: spectral theory, random Schr\"odinger operators.}

 \maketitle
\section{Introduction}
Among the most investigated and dealt with operators in the field
of mathematical physics problems are random Schr\"odinger
operators of the form
\begin{equation}
H_{\omega}=-\Delta +W_{\text{per}}+ V_{\omega}=H_{0}+ V_{\omega},
\end{equation}
where $W_{\text{per}}$ is a $\mathbb{Z}^d$-periodic function and
$V_{\omega}$ is a random potential having the Anderson form, i.e
$V_{\omega}(\cdot)=\sum_{n\in
\mathbb{Z}^d}\omega_{\gamma}f(\cdot-\gamma)$. See \cite{5,Stob},
for the physical motivations.\newline The study of the spectral
theory of operators of the form (\ref{b}) have drawn the attention
of many researchers for the importance of the related results. In
fact, it is linked to the systems evolutions for which the
Hamiltonian is described by (\ref{b}). The goal of this paper is
to discuss one of the problems that remain unsolved: the spectrum
location of $H_{\omega}$,
 precisely the spectrum infimum.
This will be carried out in the case when the single site $f$ does
not have a definite sign.\newline As the main object is to study
the location of the spectrum, let us recall the following basic
results already known on this subject and stated by Kirsch and
Martinelli \cite{5,9}:
\begin{thm}\label{bij1}
\begin{equation}
\Sigma(H_{ ,\omega})=\overline{\bigcup_{\omega_{\gamma}
\in\mathcal{P}}\Sigma\Big(-\Delta+W_{\text{per}}+
\sum_{\gamma\in\mathbb{Z}^d}
\omega_{\gamma}f(x-\gamma)\Big)}.\label{wa01}
\end{equation}
Here $\mathcal{P}$ is the set of all periodic sequences
$\{\omega_n\}_{n\in\mathbb{Z}^d}$, with an arbitrary period such
that $\omega_{n}$ is in the support of $\mu$ for all $n$ and $\Sigma(H)$ is the
spectrum of $H$.
\end{thm}
As has been said above the proof of Theorem \ref{bij1} exists in
\cite{5,9} and is based on Weyl sequences and probabilistic
arguments. Notice that this theorem reduces the determination of
the spectra of random Schr\"odinger operators for the case of
periodic Schr\"odinger operators. As it is well known \cite{Resi}
that the spectrum of periodic operators have a band structure,
this will be the case for $\Sigma(H_{\omega})$ with the
possibility to close gaps.\newline Under additional assumptions on
$f$ more is known:
\begin{thm}\label{bij2}If $f$ has a fixed sign, i.e $f\leq 0$ or $f\geq 0$ and
if $\mu$ is supported $[\omega^-,\omega^+]$, then
\begin{equation}
\Sigma(H_{ ,\omega})=\bigcup _{\omega\in [\omega^-,\omega^+]}
\Sigma\Big (-\Delta+W_{\text{per}}+
\omega\sum_{\gamma\in\mathbb{Z}^d} f(x-\gamma)\Big).\label{waa1}
\end{equation}
In particular,
\begin{equation}\inf\Sigma(H_{ ,\omega})=\left\{\begin{array}{ll}
\inf \Sigma\Big (-\Delta+W_{\text{per}}+ \omega
^+\sum_{\gamma\in\mathbb{Z}^d} f(x-\gamma)\Big)\ if\ f \leq 0,\\
\inf \Sigma\Big (-\Delta+W_{\text{per}}+ \omega ^-\sum_{\gamma\in\mathbb{Z}^d} f(x-\gamma)\Big)\ if\ f\geq 0.\\

\end{array}\right.\label{waa2}
\end{equation}
\end{thm}
Theorem \ref{bij2} is a simple consequence of Theorem \ref{bij1}.
Indeed, using (\ref{wa01}) and the fact that constant sequences
$\{\omega_{\gamma}=\omega\in \ supp \ \mu\}\in \mathcal{P}$, we
get that the r.h.s of (\ref{waa1}) is naturally contained in
$\Sigma(H_{ ,\omega})$. For the inverse inclusion, let
$(\omega_{\gamma})_{\gamma\in \mathbb{Z}^d}\in \mathcal{P}$ be
$k$-periodic and let $[a,b]$ be the $n$-th band of the $k$
periodic operator $\displaystyle
H_{\omega,k}=-\Delta+W_{\text{per}}+
\sum_{\gamma\in\mathbb{Z}^d}\omega_{\gamma}f(x-\gamma)$. Let
$[a_1,b_2]$ and $[a_2,b_2]$ be the $n$-th bands of respectively
 $-\Delta+W_{\text{per}}+ \omega^-\sum_{\gamma\in\mathbb{Z}^d}f(x-\gamma)$ and
 $-\Delta+W_{\text{per}}+ \omega^+\sum_{\gamma\in\mathbb{Z}^d}f(x-\gamma)$
  (both seen as $k$-periodic operators).
 By the min-max principle we have $a_1\leq a\leq a_2$ and $b_1\leq b\leq b_2$.
 As the bands of $\displaystyle -\Delta+W_{\text{per}}+ \omega\sum_{\gamma\in\mathbb{Z}^d}f(x-\gamma)$
 depend continuously on $\omega$, we deduce that $\Sigma\Big(-\Delta+W_{\text{per}}+
 \sum_{\gamma\in\mathbb{Z}^d}\omega f(x-\gamma)\Big)$
 is contained in the r.h.s of (\ref{waa1}). The proof of (\ref{waa1})
 is ended by taking into account the fact
 that these sets are closed.\newline
 For (\ref{waa2}) it is a simple consequence of monotonicity of the model,
 it is increasing when $f\geq 0$ and decreasing when $f\leq 0$. indeed if $f\leq 0$, and
 $\widetilde{\omega}_{\gamma}\leq \hat{\omega}_{\gamma}$ for any $\gamma\in \mathbb{Z}^d$ then
in the sense form we have
$$
-\Delta+W_{\text{per}}+\sum_{\gamma\in\mathbb{Z}^d}\hat{\omega}_{\gamma}
f(x-\gamma)\leq
-\Delta+W_{\text{per}}+\sum_{\gamma\in\mathbb{Z}^d}
\widetilde{\omega}_{\gamma}f(x-\gamma).
$$
 The situation is more complicated and different when the single site $f$ changes the
 sign, as the monotonicity property is not true in this case.
We notice that recently Lott and Stolz have conjectured
\cite{stolz} that in dimension one, the spectral minimum of random
displacement models is realized through the pair formation of the
single site.
\subsection{The Model}
Our basic object of study is the so-called \textit{Anderson
model}, a random Schr\"odinger operator of the form
\begin{equation}
H_{\lambda,\omega}=-\sum_{i=1}^d \frac{\partial^2}{\partial
x_{i}^2}+W_{\text{per}}+\lambda\sum_{\gamma\in\mathbb{Z}^d}\omega_{\gamma}f(x-\gamma),\label{b}
\end{equation}
where,
\begin{itemize}
\item $W_{\text{per}}$ is a $\mathbb{Z}^d $-periodic and bounded
function.
\item $\lambda$ is a positive parameter
\item $(\omega_{\gamma})_{\gamma\in\mathbb{Z}^d}$ is a family of
independent, identically-distributed random variables taking
values in $[\omega^{-},\omega^+]$. We denote by $\mu$ the
probability distribution.
\item For $C_{0}=[-\frac{1}{2},\frac{1}{2}]$,
the single site potential $f\in C_{0}^{\infty}(C_0)$ such that
$f\in l^{1}( L^p(\mathbb{R}^d))$, with $p=2$ if $d\leq 3$, $p>2$
if $d=4$ and $p=d/2$ if $d\geq 5$ and $f=f^++f^-$, with $f^+\geq
0,\ f^-\leq 0$ and $\ f^+\cdot \ f^-=0$.
\end{itemize}
By \cite{5,20}, we know that $H_{\lambda,\omega}$ is an
essentially self-adjoint operator on $L^2(\mathbb{R}^d)$ with the
domain $C_{0}^{\infty}(\mathbb{R}^d)$, we denote by
$H_{\lambda,\omega}$
 its self adjoint extension.\newline
It is an ergodic operator so, according to \cite{5,20}, we know
that there exist $\Sigma_{\lambda}, \Sigma _ {\lambda,pp}, \Sigma
_ {\lambda,ac}$ and $\Sigma _ {\lambda,sc}$ closed and non-random
sets of ${\Bbb{R}}$ such that $\Sigma_{\lambda}$ is the spectrum
of $H_{\lambda,\omega}$ with probability one and such that if
$\sigma_{pp}$ (respectively $\sigma_{ac}$ and $\sigma_{sc}$)
design the pure point spectrum (respectively the absolutely
continuous and singular continuous spectrum) of $H_{\lambda,
\omega}$, then $\Sigma _ {pp}=\sigma _ {pp}, \Sigma _
{\lambda,ac}=\sigma _ {\lambda,ac}$ and $\Sigma _
{\lambda,sc}=\sigma _ {\lambda,sc}$ with probability one.
\subsection{The result}
As we will see (subsection\ref{flo}) $H_{\lambda,\omega}$ can be
considered as a perturbation of some periodic operator
$H_{\lambda,\omega^-}$. Let
$\varphi_{\lambda,1}(x,\theta(\lambda))$ be the Floquet
eigenfunction associated to the first Floquet eigenvalue
$E_{1}(\lambda,\theta)$ of $H_{\lambda,\omega^-}$. Let
$(\theta_{k}(\lambda))_{1\leq k\leq m}$ be the points where
$E_1(\lambda,\theta)$ attains its minimum. We set
$$A(0)=\Big(\langle
f\varphi_{0,1}(\cdot,\theta_k(0)),\varphi_{0,1}(\cdot,\theta_{k'}(0))\rangle_{L^2(C_0)}\Big)_{1\leq
k,k'\leq m}.$$ We prove that
 \begin{thm}\label{hida}
Let $H_{\lambda,\omega}$ be the operator defined by (\ref{b}).\\
If the matrix $A(0)$ is positive-definite. Then there exists
$\lambda_{0}>0$ such that for any $\lambda\in [0,\lambda_0]$ we
have:
$$\inf(\Sigma_{\lambda})=\inf(\Sigma_{\lambda,\omega^-}).$$
If the matrix $A(0)$ is negative-definite. Then there exists
$\lambda_{0}>0$ such that for any $\lambda\in [0,\lambda_0]$ we
have:
$$\inf(\Sigma_{\lambda})=\inf(\Sigma_{\lambda,\omega^+}).$$
\end{thm}
\begin{rem}
Theorem \ref{hida} is stated for the infimum of the spectrum.
Under some additional assumptions the same result is still true
for the internal edges of the spectrum.\newline The analogous
problem for the random magnetic Schr\"odinger operator is
considered and studied in \cite{ghri}.\\ Theorem \ref{hida} can be
considered as a first step toward the physically-motivated
applications. One of them is the study of the so-called Lifshitz
tails of the integrated density of states. This could be done
under some additional assumptions on the behavior of the random
variables in the vicinity of $\omega^-$ or $\omega^+$,
\cite{12,W1,W2}. Another one is the spectral localization
\cite{Stob}.
\end{rem}
The proof of Theorem \ref{hida}, is given in section \ref{preuve}.
It is based on the reduction procedure. This powerful technique
was predicted by Klopp \cite{12} and used in several works,
\cite{ghri,W1,W2}.\newline
 As stated the proof of Theorem \ref{hida} can be divided naturally divides
 into two parts,
 we shall discuss them separately.\newline
  Indeed, if $A(0)$ is positive-definite
we will conjugate $H_{\omega}$ with $\Pi_{\lambda,0}$, the
spectral projection for $H_{\lambda,\omega^-}$ on the first band.
Then we prove that
$$\Pi_{\lambda,0}H_{\lambda,\omega}\Pi_{\lambda,0}\geq
E_{\lambda,\omega^-}\Pi_{\lambda,0}.$$ Here $E_{\lambda,\omega^-}$
is the bottom of the spectrum of the periodic operator
$H_{\lambda,\omega^-}$. \newline If $A(0)$ is negative-definite we
will conjugate $H_{\omega}$ with $\Pi_{\lambda,0}$, the spectral
projection for $H_{\lambda,\omega^+}$ on the first band. Then we
prove $$\Pi_{\lambda,0}H_{\lambda,\omega}\Pi_{\lambda,0}\geq
E_{\lambda,\omega^+}\Pi_{\lambda,0}.$$ Here $E_{\lambda,\omega^+}$
is the bottom of the spectrum of the periodic operator
$H_{\lambda,\omega^+}$.
 \section{Preliminary}
 Let us consider the following periodic operator
 \begin{equation}
 H_{\lambda,\omega^-}=-\Delta +W_{\text{per}}+\lambda\sum_{\gamma\in\mathbb{Z}^d}\omega^-f(\cdot-\gamma)
.\end{equation}
 For this, it is convenient to consider $H_{\lambda,\omega}$
 as a perturbation of $H_{\lambda,\omega^-}$. Indeed,
 we have:
 $$
 H_{\lambda,\omega}=H_{\lambda,\omega^-}+\lambda\sum_{\gamma\in\mathbb{Z}^d}(\omega_{\gamma}-\omega^{-})f(\cdot-\gamma).
 $$
 For this for $\gamma\in\mathbb{Z}^d$, we set
$\widetilde{\omega}_{\gamma}=\omega_{\gamma}-\omega^-$ and
$V_{\widetilde{\omega}}(\cdot)=\sum_{\gamma\in\mathbb{Z}^d}\widetilde{\omega}_{\gamma}f(\cdot-\gamma)$.
  We notice that according to the definition of
  $(\omega_{\gamma})_{\gamma\in\mathbb{Z}^d}$ we get that, $(\widetilde{\omega}_{\gamma})_{\gamma\in\mathbb{Z}^d}$ is a
family of random positive and bounded variables. \subsection{Some
Floquet Theory}\label{flo}
 For $\gamma\in\mathbb{Z}^d$, we denote by $\tau_{\gamma}$ the translation by $\gamma$ operator
i.e $(\tau_{\gamma}\varphi)(x)=\varphi(x-\gamma)$. We have, for
any $\gamma\in \mathbb{Z}^d$
$$
\tau_{\gamma}H_{\lambda,\omega^-}\tau_{\gamma}^*=
\tau_{\gamma}^*H_{\lambda,\omega^-}\tau_{\gamma}=H_{\lambda,\omega^-}.
 $$
Then the so called, Floquet Theory, can be used to study
$H_{\lambda,\omega^-}$.
 For this, we review some standard
facts from the Floquet theory for periodic operators. Basic
references for this material are \cite{15,Resi,23}. Let
$\Bbb{T}^*=\Bbb{R}^d/(2\pi{\Bbb{Z}^d})$. We define
${{\mathcal{H}}}$ by
$$
{{\mathcal{H}}}=\{u(x,\theta)\in L^{2}_{loc}({\Bbb{R}}^d)\otimes L^2({\Bbb{T}%
}^*); \forall (x,\theta,\gamma)\in \Bbb{R}^d\times \Bbb{T}^*\times \Bbb{Z}%
^d;\ u(x+\gamma,\theta)=e^{i\gamma\theta}u(x,\theta)\}.
$$
$\mathcal{H}$ is equipped with the norm
$$
\frac{1}{\text{vol}(\mathbb{T}^*)}\int_{\mathbb{T}^*}\|u(x,\theta)\|_{L^{2}(C_{0})}^2d\theta.
$$
 For $\theta\in \mathbb{R}^d$ and $u\in\mathcal{S}(\mathbb{R}^d)$;
the Schwartz space of rapidly decreasing functions we define
$$
(Uu)(x,\theta)=\sum_{\gamma\in \mathbb{Z}^d}e^{i\gamma\cdot
\theta}u(x-\gamma).
$$
$U$ can be extended as a unitary isometry from
$L^{2}(\mathbb{R}^d)$ to $\mathcal{H}$. Its inverse is given by
the formula,
$$
\text{for} \ u\in \mathcal{H},\ \
(U^*u)(x)=\frac{1}{\text{vol}(\mathbb{T}^*)}\int_{\mathbb{T}^*}u(x,\theta)d\theta.
$$
$U$ is a unitary isometry from $L^{2}({\Bbb{R}}^d)$ to ${{\mathcal{%
H}}}$ and $H_{\lambda,\omega}$ admits the Floquet decomposition
\cite{15,23}
$$
UH_{\lambda,\omega^-}U^*=\frac{1}{\text{vol}(\mathbb{T}^*)}\int_{{\Bbb{T}^*}}^{\oplus}H_{\lambda,\omega-}(\theta)d\theta.
$$
Here $H_{\lambda,\omega^-}(\theta)$ is the operator $H_{\lambda,\omega^-}$ acting on ${{\mathcal{H}}}_{\theta}$,
defined by
$$
{\mathcal{H}}_{\theta}=\{u\in L^2_{loc}(\Bbb{R}^d); \forall \gamma\in {\Bbb{Z%
}^d}, u(x+\gamma)= e^{i\gamma\theta}u(x)\}.
$$
As $H_{\lambda,\omega^-}$ is elliptic, we know that,
$H_{\lambda,\omega^-}(\theta)$ has a compact resolvent; hence its
spectrum is discrete \cite{Resi}. We denote its eigenvalues,
called Floquet eigenvalues of $H_{\lambda,\omega^-}$, by
$$
E_{1}(\lambda,\theta)\leq E_{2}(\lambda,\theta)\leq \cdot
\cdot\cdot\leq E_{n}(\lambda,\theta)\leq \cdot\cdot \cdot.
$$
The corresponding Floquet eigenfunctions are denoted by $(\varphi_{\lambda,j}(x,\cdot))_{j\in {%
\Bbb{N}}^*}$. The functions $(\theta\to
E_{n}(\lambda,\theta))_{n\in {\Bbb{N}}^*}$ are
Lipshitz-continuous, and we have
$$
E_n(\lambda,\theta)\to +\infty \ \ \mathrm{as}\ n\to +\infty\ \ \mathrm{uniformly \
in}\ \ \theta.
$$
The spectrum $\Sigma_{\lambda,\omega^-}$ of $H_{\lambda,\omega^-}$ is made of bands (i.e $%
\displaystyle \Sigma_{\lambda,\omega^-}=\cup_{n\in
\Bbb{N}^*}E_{n}(\lambda,\Bbb{T}^*). $)\newline Let us note by
$E_{\lambda,\omega^-}$ the bottom of the spectrum of
$\Sigma_{\lambda,\omega^-}$,
 i.e $E_{\lambda,\omega^-}=\inf_{\theta\in\mathbb{T}^d}E_{1}(\lambda,\theta)$.\newline
 It is a well-known fact that, in any dimension the bottom (the first band) of the spectrum of a periodic Schr\"odinger
 operators is given by a simple Floquet eigenvalue and that the minimum of this Floquet eigenvalue is non-degenerate and quadratic. More
 precisely let $\theta(\lambda)$ be an element of
 $$
 Z_{\lambda}=\{\theta\in\mathbb{T}^*;\ E_{1}(\lambda,\theta)=E_{\lambda,\omega^-}\}.
 $$
 Then there exist $C>0$ and $\delta>0$ such that
 $$
 \forall \ |\theta-\theta(\lambda)|\leq \delta,\ \ \ \ \frac{1}{C} |\theta-\theta(\lambda)|^2\leq E_{1}(\lambda,\theta)-E_{\lambda,\omega^-}\leq
 C|\theta-\theta(\lambda)|^2.
 $$
Hence, the points where $E_1(\lambda,\theta)$ reaches
$E_{\lambda,\omega^-}$ are isolated and as $\mathbb{T}^*$ is
compact, one concludes that $Z_{\lambda}$ contains only finitely
many of elements. Let $m$ be the cardinal of $Z_{\lambda}$ and let
us denote them by $(\theta_{k}(\lambda))_{1\leq k\leq m}$. One can
check that $\theta_k(\lambda)$ depends continuously on $\lambda$.
For the sake of brevity, we use the notation
$\theta_k=\theta_k(\lambda)$.\newline For $1\leq k\leq m$ and
$\theta\in\mathbb{T}^*$, we set
\begin{equation}
\zeta_{k,\lambda}(\theta)=\sum_{1\leq i\leq
d}(\theta_i-\theta_{k,i})^2.\label{joum}
\end{equation}
We notice that for any $1\leq k\leq m$,
$$
\theta \mapsto \varphi_{\lambda,1}(x,\theta),
$$
is analytic in a neighborhood of $\theta_{k}$.
\subsection{Wannier basis}
We recall concepts used in \cite{12}. Let ${{\mathcal{E}}}\subset
L^2({\Bbb{R}} ^d)$ be a closed subspace invariant by the
${\Bbb{Z}^d}$-translations, i.e.
such that $\Pi^{{\mathcal{E}}} $, the orthogonal projection on ${\mathcal{E}%
}$, satisfies
$$
\forall {\gamma} \in {\Bbb{Z}^d}, \
\Pi^{{\mathcal{E}}}=\tau_{-\gamma }\Pi
^{{\mathcal{E}}}\tau_{\gamma}.
$$
Following the computations done in section 1.2 of \cite{12}, we
see that there exists an orthonormal system of vectors
$(\widetilde{\varphi}_{j, 0})_{j\in
N}$ such that for $\widetilde{\varphi}_{j, {\gamma} }=\tau_{\gamma}(\widetilde{\varphi}%
_{j, 0})$; $(\widetilde{\varphi}_{j, {\gamma}})_{(j\in
N;{\gamma}\in {\Bbb{Z}^d})}
$ is an orthonormal basis of ${\mathcal{E}}$. Such system is called \textbf{%
Wannier basis} of $\mathcal{\ E}$. The vectors
$(\widetilde{\varphi}_{n, 0})_{n\in N}$ are called \textbf{Wannier
generators} of ${\mathcal{E}}$. \newline Let ${{\mathcal{E}}}
\subset L^2({\Bbb{R}} ^d)$ be a space which is
translation-invariant. ${\mathcal{E}}$ is said to be of finite energy for $%
H_{\lambda,\omega}$ if $\Pi^{{\mathcal{E}}} H_{\lambda,\omega}
\Pi^{\mathcal{\ E}}$ is a bounded operator. In this case,
${\mathcal{E}}$ admits a finite set of Wannier generators.\newline
Let $\Pi_{\lambda,0}(\theta)$ (respectively
$\Pi_{\lambda,+}(\theta)$) be the orthogonal projection in
$\mathcal{H}_{\theta}$ on the vector space generated by
$\varphi_{\lambda,1}(\cdot,\theta)$ (respectively by
$(\varphi_{\lambda,j}(\cdot,\theta))_{j\geq 2}$ ). These
projections are two-by-two orthogonal and their sum is the
identity for all $\theta \in {\Bbb{T}^*}$. One defines
$$
\Pi_{\lambda,\alpha} = U^{-1}\Big
(\int_{\Bbb{T}^*}\Pi_{\lambda,\alpha}(\theta)d\theta \Big) U:
L^2({\Bbb{R}}^d)\to L^2({\Bbb{R}}^d),
$$
where $\alpha \in \{ 0, +\}$. $\Pi_{\lambda,\alpha}$ is an
orthogonal projection
on $L^2(\Bbb{R}^d) $ and for all ${\gamma} \in {\Bbb{Z}^d}$, we have $%
\tau^*_{\gamma
}\Pi_{\lambda,\alpha}\tau_{\gamma}=\Pi_{\lambda,\alpha}.$ \newline
For $\alpha \in \{0, +\}$, we set ${\mathcal{%
E}}_{\lambda,\alpha} = \Pi_{\lambda,\alpha}(L^2({\Bbb{R}} ^d)).$
These spaces are translation-invariant. Moreover
${\mathcal{E}}_{\lambda,0}$ is of finished energies for
$H_{\lambda,\omega^-}$. The reduction procedure consists in
decomposing the operator $H_{\lambda,\omega^-}$ according to
various translation-invariants subspaces. The random operators
thus obtained are reference operators.
\section{The proof of Theorem \ref{hida}\label{preuve}}
As we have indicated, our aim in this section is to prove Theorem
\ref{hida}, but first, let us introduce some notations and useful
lemma.\newline For $u\in L^{2}(\mathbb{T}^*)$, let
\begin{equation}
\mathcal{T}_{\varphi_{\lambda,1}}(u)=U^*(u\varphi_{\lambda,1}(x,\theta))=\int_{\mathbb{T}^*}u(\theta)\varphi_{\lambda,1}
(x,\theta)d\theta.\label{per1}
\end{equation}
So, ${\mathcal{T}}_{\varphi_{\lambda,1}}$ define a unitary
transformation from $L^2(\mathbb{T}^*)$ to
${\mathcal{E}}_{\lambda,0}$ and for $v\in \mathcal{E}_{\lambda,0}$
we have
$$
{\mathcal{T}}_{\varphi_{\lambda,1}}^*(v)=\langle
(Uv)(\cdot,\theta),\varphi_{\lambda,1}(\cdot,\theta)\rangle.
$$
For $1\leq k\leq m$ and $(x,\theta)\in \mathbb{R}^d\times
\mathbb{T}^*$, let
$$
\widetilde{\varphi}_{\lambda,1,k}(x,\theta)=\varphi_{\lambda,1}(\theta_{k}(\lambda),x)e^{(\theta-\theta_k)\cdot
x}.
$$
We set,
$$\delta\varphi_{\lambda,1,k}(x,\theta)=\frac{1}{\sqrt{\zeta_{k,\lambda}(\theta)}}\Big(\varphi_{\lambda,1}-\varphi_{\lambda,1,k}(x,\theta)\Big).$$
By this, for any $u\in\L^2(\mathbb{T}^*)$, we have
\begin{equation}
\mathcal{T}_{\varphi_{\lambda,1}}(u)=\mathcal{T}_{\widetilde{\varphi}_{\lambda,1,k}}(u)
+\mathcal{T}_{\delta\varphi_{\lambda,1,k}}(\sqrt{\zeta_{k,\lambda}}u).\label{wakt1}
\end{equation}For $\displaystyle v\in {\mathcal{H}}_{\theta}^2=\{v\in H^{2}_{\text{loc}}(\mathbb{R}^d);
 v(\cdot-\gamma)=e^{-i\gamma\cdot \theta}u(\cdot)\}$ one defines the following
norms:
$$\sup_{\theta\in
\mathbb{T}^*}\{\|v(\cdot,\theta)\|_{L^{2}(C_0)}^2\}=\|v\|_{1,\infty}^2.
$$
and
$$
\sup_{\theta\in\ {\Bbb{T}^*}}\Big(\|H_{{\lambda,\omega^-}%
}(\theta)u(\cdot,\theta)\|_{L^2(C_{0})}^2+ ||v(\cdot,\theta)||^2_{L^2(C_{0})}%
\Big)=||v||_{H_{\lambda,\omega^-},\infty}.
$$
\begin{rem}
The functions $\displaystyle\widetilde{\varphi}_{\lambda,1,k}$ and
$\displaystyle \delta\varphi_{\lambda,1,k}$ are well defined and
$$\|\widetilde{\varphi}_{\lambda,1,k}\|_{1,\infty},\
\|\widetilde{\varphi}_{\lambda,1,k}\|_{H_{\lambda,\omega^-},\infty}\
\|\delta\widetilde{\varphi}_{\lambda,1,k}\|_{1,\infty}, \
{\text{and}}\
\|\delta\widetilde{\varphi}_{\lambda,1,k}\|_{H_{\lambda,\omega^-},\infty};$$
are finished (See \cite{12}).
\end{rem}
The following Lemma is of use. It will be proven at the end of
this section.
\begin{lem}\label{sousou1}
For $\theta_k,\theta_{k'}\in \mathbb{T}^*$ and $\varphi\in
L^2(\mathbb{T}^*,\mathcal{H}^2)$ let $
\displaystyle\varphi_{k}=e^{i(\theta-\theta_k)x}\varphi(x,\theta_k)$,
$\displaystyle
\varphi_{k'}=e^{i(\theta-\theta_{k'})x}\varphi(x,\theta_{k'})$ and
$\displaystyle
a_{\varphi,k,k'}(x)=f(x)\varphi(x,\theta_k)\varphi(x,\theta_{k'})$.
If $\displaystyle\|\varphi\|_{1,\infty}<\infty$ (resp.
$\displaystyle\|\varphi\|_{H_{\lambda,\omega^-},\infty}<\infty$ )
then $\displaystyle\mathcal{T_{\varphi}}\in
\mathcal{L}(L^{2}(\mathbb{T}^*),L^{2}(\mathbb{R}^d))$ (resp.
$\displaystyle V_{\widetilde{\omega}}\cdot\mathcal{T_{\varphi}}\in
\mathcal{L}(L^{2}(\mathbb{T}^*),L^{2}(\mathbb{R}^d))$) and there
exist $C,\beta>0$ such that for all $\displaystyle u,v \in
L^{2}(\mathbb{T}^*)$, we have
\begin{multline}
\Big|\langle V_{\widetilde{\omega}}{\mathcal{T}}_{\varphi_{k}}
(u),{\mathcal{T}}_{\varphi_{k'}}(v)\rangle-\Big(\int_{C_0}a_{\varphi,k,k'}(x)dx\Big)\cdot
\sum_{\gamma\in\mathbb{Z}^d}\widetilde{\omega}_{\gamma}\hat{u}(\gamma)\overline{\hat{v}(\gamma)}\Big|\\
\leq
C\beta\sum_{\gamma\in\mathbb{Z}^d}\widetilde{\omega}_{\gamma}(|\hat{u}(\gamma)|^2+|\hat{v}(\gamma)|^2)+C(1+1/\beta)(\langle
\zeta_{k,\lambda}u,u\rangle_{L^{2}(\mathbb{T}^*)}+\langle
\zeta_{k',\lambda}v,v\rangle_{L^{2}(\mathbb{T}^*)}).\label{eq03}
\end{multline}
\begin{equation}
\|V_{\widetilde{\omega}}{\mathcal{T}}_{\varphi_{k}}(u)\|^{2}_{L^{2}(\mathbb{R}^d)}\leq
C\Big(\sum_{\gamma\in\mathbb{Z}^d}\widetilde{\omega}_{\gamma}|\hat{u}(\gamma)|^2+
\langle\zeta_{k,\lambda}u,u\rangle_{L^{2}(\mathbb{T}^*)}\Big).\label{eq04}
\end{equation}
\end{lem}
\subsection{If $A(0)$ is positive-definite\label{tak}}
 We set
$$H_{\lambda,\omega}^0=\Pi_{\lambda,0}H_{\lambda,\widetilde{\omega}}\Pi_{\lambda,0}=
\Pi_{\lambda,0}H_{\lambda,\omega^-}\Pi_{\lambda,0}+\Pi_{\lambda,0}V_{\widetilde{\omega}}\Pi_{\lambda,0}.$$

\begin{thm}\label{rhou}
Assume that the matrix $A(0)$ is positive-definite. Then there
exists $\lambda_{0}>0$ such that for any $\lambda\in
[0,\lambda_0]$ we have:
$$\Pi_{\lambda,0}\Big(H_{\lambda,\omega}-E_{\lambda,\omega^-}\Big)\Pi_{\lambda,0}\
\ \ \text{is\ a \ positive\ operator.}
$$
\end{thm}
The proof of Theorem \ref{rhou} is the object of the following
section. \subsubsection{The proof of Theorem \ref{rhou}} Using
(\ref{per1}), we get that $H_{\lambda,\omega}^{0} $ is unitarly
equivalent to the operator
$$h_{\lambda,{\omega}}^0=\mathcal{T}^*_{\varphi_{\lambda,1}}
H_{\lambda,{\omega}}^0{\mathcal{T}}_{\varphi_{\lambda,1}},$$
acting on $L^2(\mathbb{T}^*)$ and written as
$$
h_{\lambda,{\omega}}^0=h_{\lambda,\omega^-}^0+\lambda
V_{\lambda,\widetilde{\omega}}^0.
$$
With $h_{\lambda,\omega^-}^0$ is the multiplication operator by
$E_{1}(\lambda,\theta)$ and $V_{\lambda,\widetilde{\omega}}^0$ is
an integral operator with the kernel
$$
V_{\lambda,\widetilde{\omega}}(\theta,\theta ')=\langle
V_{\widetilde{\omega}}\varphi_{\lambda,1}(\cdot,\theta),\varphi_{\lambda,1}
(\cdot,\theta')\rangle.
$$
Let $V_{k}$ be a neighborhood of $\theta_{k}$, such that if
$\theta_{k'}\in Z$ and $k\neq k'$ then $\theta\notin
\overline{V_{k}}$ and $V_k\cap V_{k'}=\emptyset$. As
$\mathbb{T}^*$ is compact, one can cover it by $(V_{k})_{1\leq
k\leq m}$ (i.e $\displaystyle\cup_{1\leq k\leq
m}V_{k}=\mathbb{T}^*$). For $1\leq k\leq m$ let $\chi_k$ be the
characteristic function of $V_k$.\\ For simplicity for $u\in
L^2(\mathbb{T}^*)$, we will denote $\chi_ku$ as $u_k$ in the
following. We consider $u$ as a system of $m$ columns denoted by
$(u_{ k})_{1\leq k\leq m}$. We endow $L^2({\Bbb{T}^*})\otimes
{\Bbb{C}}^m $ with the scalar product generating the following
Euclidean norm:
$$
\| u\| _{L^2({\Bbb{T}^*})\otimes {\Bbb{C}}%
^m}^2=\sum_{k=1} \| u_{k}\| _{L^2({\Bbb{T}^*})}^2.
$$
\subsubsection{The lower bound of $h_{\lambda,\omega^-}^0$}
\begin{prop}\label{wakt2}
There exists $C>0$ such that for any $u\in L^2(\mathbb{T}^*)$, we
have
\begin{equation}
\langle h_{\lambda,\omega^-}^0u,u\rangle\geq \sum_{1\leq k\leq
m}E_{1}(\lambda,\theta_{k}(\lambda))\|u_k\|^2_{L^{2}(\mathbb{T}^*)}+\frac{1}{C}\sum_{1\leq
k\leq
m}\langle\zeta_{k,\lambda}u_k,u_k\rangle_{L^{2}(\mathbb{T}^*)}
.\label{wa1}
\end{equation}
\end{prop}
\noindent {\bf{Proof:}} For $u\in L^2(\mathbb{T}^*)$, one computes
\begin{equation*}\langle
h_{\lambda,\omega^-}^0u,u\rangle=\int_{\mathbb{T}^*}E_1(\lambda,\theta)|u(\theta)|^2d\theta
=\sum_{1\leq k \leq
m}\int_{\mathbb{T}^*}E_{1}(\lambda,\theta)\chi_{k}(\theta)|u(\theta)|^2d\theta,
\end{equation*}
As for any $\theta $ in $V_k$ the support of $\chi_{k}$, there
exists $C>0$ such that we have
$$
E_{1}(\lambda,\theta_{k}(\lambda))+\frac{1}{C}\zeta_{k,\lambda}(\theta)\leq
E_{1}(\lambda,\theta),
$$
we get the result.\hfill$\Box$ \subsubsection{The lower bound of
$V_{\lambda,\widetilde{\omega}}^0$}
\begin{prop}\label{permi}
There exists $C_1,C_2>0$ and $\lambda_0$ such that for all
$\lambda\in [0,\lambda_0]$ and $u\in L^2(\mathbb{T}^*)$ we have
\begin{equation}
\langle V_{\lambda,\widetilde{\omega}}^0u,u\rangle \geq
C_1\sum_{1\leq k\leq m,\gamma \in
\mathbb{Z}^d}\widetilde{\omega}_{\gamma}|(\hat{u_k})(\gamma)|^2-C_2\sum_{1\leq
k\leq
m}\langle\zeta_{k,\lambda}u_k,u_k\rangle_{L^{2}(\mathbb{T}^*)}.
\end{equation}
\end{prop}
{\bf{The proof of Theorem \ref{rhou}}} \\ \noindent Let us notice
that, by combining the results of Proposition \ref{wakt2} and
\ref{permi}, one gets that there exists $\lambda_0>0$ such that
for any $\lambda\in[0,\lambda_0]$ and for any $u\in
L^{2}(\mathbb{T}^*$ we have,
\begin{eqnarray*}
\langle h_{\lambda,\omega}^0u,u\rangle&\geq& \sum_{1\leq k\leq
m}E_{1}(\lambda,\theta_{k}(\lambda))\|u_k\|^{2}_{L^2(\mathbb{T}^*)}+\frac{1}{C}\Big(\sum_{1\leq
k\leq
m}\langle\zeta_{k,\lambda}u_k,u_k\rangle_{L^2(\mathbb{T}^*)}+\\
& &\lambda\sum_{1\leq k\leq
m}\sum_{\gamma\in\mathbb{Z}^d}\widetilde{\omega}_{\gamma}|\hat{u_k}(\gamma)|^2\Big)\\
&\geq&E_{\lambda,\omega^-}\|u\|_{L^{2}(\mathbb{T}^*)}^2+\frac{1}{C}\Big(\sum_{1\leq
k\leq
m}\langle\zeta_{k,\lambda}u_k,u_k\rangle_{L^2(\mathbb{T}^*)}+\\
& &\lambda\sum_{1\leq k\leq
m}\sum_{\gamma\in\mathbb{Z}^d}\widetilde{\omega}_{\gamma}|\hat{u_k}(\gamma)|^2\Big).
\end{eqnarray*}
This gives that,
\begin{multline}
\langle (
h_{\lambda,\widetilde{\omega}}^0-E_{\lambda,\omega^-})u,u\rangle
\geq\\\frac{1}{C}\Big(\sum_{1\leq k\leq
m}\langle\zeta_{k,\lambda}u_k,u_k\rangle_{L^2(\mathbb{T}^*)}+\lambda\sum_{1\leq
k\leq
m}\sum_{\gamma\in\mathbb{Z}^d}\widetilde{\omega}_{\gamma}|\hat{u_k}(\gamma)|^2\Big).
\end{multline}
This ends the proof of Theorem \ref{rhou}.\hfill$\Box$\newline
\begin{rem}
We notice that even if we know that the bottom of the spectrum of
$H_{\omega}$ coincides with the bottom of the spectrum of
$H_{\lambda,\omega^-}$ we cannot consider
$\Pi_{\lambda,0}\Big(H_{\lambda,\omega}-H_{\lambda,\omega^-}\Big)\Pi_{\lambda,0}$
as a positive operator.
\end{rem}
\noindent {\bf{The proof of proposition \ref{permi}:}} \\
\noindent Let us start by expanding $\langle
V_{\lambda,\widetilde{\omega}}^0u,u\rangle$,
\begin{eqnarray*} \langle V_{\lambda,\widetilde{\omega}}^0u,u\rangle
&=&\sum_{1\leq k,k'\leq m}\langle
V_{\lambda,\widetilde{\omega}}^0{\mathcal{T}}_{\widetilde{\varphi}_{\lambda,1,k}}(u_k),{\mathcal{T}}_{\widetilde{\varphi}_{\lambda,1,k'}}(u_{k'})\rangle_{L^2(\mathbb{R}^d)}\\
&+&\sum_{1\leq k,k'\leq m}\langle
V_{\lambda,\widetilde{\omega}}^0{\mathcal{T}}_{\delta\widetilde{\varphi}_{\lambda,1,k}}(\sqrt{\zeta_{k,\lambda}}u_k),
{\mathcal{T}}_{\delta\widetilde{\varphi}_{\lambda,1,k'}}(\sqrt{\zeta_{k',\lambda}}u_{k'})\rangle_{L^2(\mathbb{R}^d)}\\
&+&2\sum_{1\leq k,k'\leq m}\Re\Big(\langle
V_{\lambda,\widetilde{\omega}}^0
{\mathcal{T}}_{\widetilde{\varphi}_{\lambda,1,k}}(u_k),
{\mathcal{T}}_{\delta\widetilde{\varphi}_{\lambda,1,k'}}(\sqrt{\zeta_{k',\lambda}}
u_{k'})\rangle_{L^2(\mathbb{R}^d)}\Big).
\end{eqnarray*} We start by estimating the three sums of the last equation.\newline
For the second sum, using Cauchy Schwartz inequality and Lemma
\ref{sousou1}, we get that for any $1\leq k,k'\leq m$, there
exists $C>0$ such that we have, \begin{multline*}\Big |\langle
V_{\lambda,\widetilde{\omega}}^0
{\mathcal{T}}_{\delta\widetilde{\varphi}_{\lambda,1,k}}(\sqrt{\zeta_{k,\lambda}}u_k),
{\mathcal{T}}_{\delta\widetilde{\varphi}_{\lambda,1,k'}}(\sqrt{\zeta_{k',\lambda}}u_{k'})
\rangle_{L^2(\mathbb{R}^d)} \Big|\\ \leq
 \ \ \frac{1}{2}\Big(\|V_{\lambda,\widetilde{\omega}}^0
{\mathcal{T}}_{\delta\widetilde{\varphi}_{\lambda,1,k}}(\sqrt{\zeta_{k,\lambda}}u_k)\|^2_{L^2(\mathbb{R}^d)}
+\|{\mathcal{T}}_{\delta\widetilde{\varphi}_{\lambda,1,k'}}(\sqrt{\zeta_{k',\lambda}}u_{k'})
\|^2_{L^2(\mathbb{R}^d)} \Big )\end{multline*}$$ \leq C\cdot \Big(
\langle \zeta_{k,\lambda}u_k,u_k\rangle_{L^{2}(\mathbb{T}^*)} +
\langle
\zeta_{k',\lambda}u_{k'},u_{k'}\rangle_{L^{2}(\mathbb{T}^*)}\Big).$$
So there exists $C>0$ such that we have
\begin{multline}\sum_{1\leq k,k'\leq m}\Big |\langle
V_{\lambda,\widetilde{\omega}}^0
{\mathcal{T}}_{\delta\widetilde{\varphi}_{\lambda,1,k}}(\sqrt{\zeta_{k,\lambda}}u_k),
{\mathcal{T}}_{\delta\widetilde{\varphi}_{\lambda,1,k'}}(\sqrt{\zeta_{k',\lambda}}u_{k'})
\rangle_{L^2(\mathbb{R}^d)} \Big|\\ \leq C\sum_{1\leq k\leq
m}\langle
\zeta_{k,\lambda}u_k,u_k\rangle.\label{rap1}\end{multline} For the
third sum, using the Cauchy-Schwartz inequality once more, we get
that for $1\leq k\leq m$, there exists $\beta>0$ such that we have
\begin{multline}\Big | \langle
V_{\lambda,\widetilde{\omega}}^0{\mathcal{T}}_{\widetilde{\varphi}_{\lambda,1,k}}
(u_k),{\mathcal{T}}_{\delta\widetilde{\varphi}_{\lambda,1,k'}}(\sqrt{\zeta_{k',\lambda}}u_{k'})\rangle_{L^2(\mathbb{R}^d)}
\Big|\\ \leq \beta
\|V_{\widetilde{\omega}}^0\mathcal{T}_{\widetilde{\varphi}_{\lambda,1,k}}(u_k)\|^{2}_{L^{2}(\mathbb{R}^d)}
+1/(4\beta)\|\mathcal{T}_{\delta\widetilde{\varphi}_{\lambda,1,k'}}(\sqrt{\zeta_{k,\lambda}}u_{k'})\|^{2}_{L^{2}(\mathbb{R}^d)}
.\label{rap2}\end{multline} Using equation (\ref{eq04}), one gets
that there exist $\widetilde{C}_1,\widetilde{C}_2>0$ such that
\begin{multline}\sum_{1\leq k,k'\leq m}\Big | \langle
V_{\lambda,\widetilde{\omega}}^0{\mathcal{T}}_{\widetilde{\varphi}_{\lambda,1,k}}
(u_k),{\mathcal{T}}_{\delta\widetilde{\varphi}_{\lambda,1,k'}}(\sqrt{\zeta_{k',\lambda}}u_{k'})\rangle_{L^2(\mathbb{R}^d)}
\Big|\\ \leq \widetilde{C}_1\beta\sum_{1\leq k\leq
m}\sum_{\gamma\in\mathbb{Z}^d}\widetilde{\omega}_{\gamma}|\hat{u_k}(\gamma)|^2+\widetilde{C}_2(\beta+1/\beta)\sum_{1\leq
k\leq m}\langle
\zeta_{\lambda,k}u_k,u_k\rangle_{L^{2}(\mathbb{T}^*)}.\label{rap3}\end{multline}
 For the first sum, using (\ref{eq03}), we get that there exist
 $C'_1,C'_2>0$ such that
\begin{multline}
\Big| \sum_{1\leq k,k'\leq m}\langle
V_{\widetilde{\omega}}^0\mathcal{T}_{\widetilde{\varphi}_{\lambda,1,k}}(u_k),
\mathcal{T}_{\widetilde{\varphi}_{\lambda,1,k'}}(u_{k'})\rangle_{L^{2}(\mathbb{R}^d)}
\\ -\sum_{1\leq k,k'\leq m}\sum_{\gamma\in \mathbb{Z}^d}\widetilde{\omega}_{\gamma}
\big(\int_{C_0}a_{\varphi_{\lambda,1},k,k'}(x)dx\big)(\hat{u_k})(\gamma)
\overline{(\hat{u_{k'}})(\gamma)}\Big| \\ \leq
C'_1\beta\sum_{1\leq k\leq
m}\sum_{\gamma\in\mathbb{Z}^d}\widetilde{\omega}_{\gamma}|\hat{u_k}(\gamma)|^2+C'_2(1+1/\beta)\sum_{1\leq
k\leq m}\langle
\zeta_{\lambda,k}u_k,u_k\rangle_{L^{2}(\mathbb{T}^*)}.
\label{yas2}\end{multline} Now equations (\ref{rap1}),
(\ref{rap3}) and (\ref{yas2}) give that there exist $K_1,K_2>0$
such that
\begin{multline}
\Big|\langle V_{\lambda,\widetilde{\omega}}^0u,u\rangle
-\sum_{1\leq k,k'\leq m}\sum_{\gamma\in
\mathbb{Z}^d}\widetilde{\omega}_{\gamma}
\Big(\int_{C_0}a_{\varphi_{\lambda,1},k,k'}(x)dx\Big)(\hat{u_k})(\gamma)
\overline{(\hat{u_{k'}})(\gamma)}\Big|\leq \\
K_1\beta\sum_{1\leq k\leq m}\sum_{\gamma\in
\mathbb{Z}^d}\widetilde{\omega}_{\gamma}|(\hat{u_k})(\gamma)|^2+K_2(1+1/\beta)\sum_{1\leq
k\leq m}\langle
\zeta_{k,\lambda}u_k,u_k\rangle_{L^{2}(\mathbb{T}^*)}.\label{wakt5}
\end{multline}
Now, if the matrix
$$
A(0)=\Big(\int_{C_0}a_{\varphi_{0,1},k,k'}(x)dx\Big)_{1\leq
k,k'\leq m},
$$
is positive-definite, one gets that $\inf \sigma(A(0))=C>0$
satisfies
$$
CI_{m}\leq A.
$$
Let $A(\lambda)$ be the matrix,
$$\Big(\int_{C_0}a_{\varphi_{\lambda,1},k,k'}(x)dx\Big)_{1\leq k,k'\leq m}. $$
Notice that for any $1\leq k,k'\leq m$, the functions
$$
f_{k,k'}:\lambda \to
\int_{C_0}a_{\varphi_{\lambda,1},\theta_{k}(\lambda),\theta_{k'}(\lambda)}(x)dx,$$
are continuous in $\lambda$. So there exists $\lambda_0>0$ such
that for any $\lambda\in [0,\lambda_0]$
$$
\frac{C}{2}I_m \leq A(\lambda).
$$
This gives that for any $u\in L^{2}(\mathbb{T}^*)$
\begin{multline}\frac{C}{2}\sum_{1\leq k\leq
m}\sum_{\gamma\in \mathbb{Z}^d}\widetilde{\omega}_{\gamma}|(\hat{u_k})(\gamma)|^2\\
\leq \sum_{1\leq k,k'\leq m}\sum_{\gamma\in
\mathbb{Z}^d}\widetilde{\omega}_{\gamma}\Big(\int_{C_0}
a_{\varphi_{\lambda,1},{k},{k'}}(x)dx\Big)(\hat{u_k})(\gamma)\overline{(\hat{u_{k'}})(\gamma)}.
\label{yas1}\end{multline} Now using the expanation of $\langle
V_{\lambda,\widetilde{\omega}}^0u,u\rangle $ and equation
(\ref{wakt5}), we get that there exist $K_1,K_2>0$ and $\beta>0$
such that
\begin{multline}\langle V_{\lambda,\widetilde{\omega}}^0u,u\rangle
\geq\\
(\frac{C_1'}{2}-K_1\beta)\sum_{1\leq k\leq m;\gamma\in
\mathbb{Z}^d}\widetilde{\omega}_{\gamma}|(\hat{u_k})(\gamma)|^2-K_2(1+1/\beta)\sum_{1\leq
k\leq
m}\langle\zeta_{k,\lambda}u_k,u_k\rangle_{L^{2}(\mathbb{T}^*)}.
\end{multline}So, for $\beta>0$, well chosen we get that there
exist constants $C_1,C_2>0$ such that
\begin{equation}\langle V_{\lambda,\widetilde{\omega}}^0u,u\rangle
\geq \frac{C_1}{3}\sum_{1\leq k\leq m;\gamma\in
\mathbb{Z}^d}\widetilde{\omega}_{\gamma}|(\hat{u_k})(\gamma)|^2-C_2\sum_{1\leq
k\leq
m}\langle\zeta_{k,\lambda}u_k,u_k\rangle_{L^{2}(\mathbb{T}^*)}.
\end{equation}
This ends the proof of Proposition \ref{permi}.\hfill $\Box$
\subsection{If $A(0)$ is negative-definite}Let
$H_{\lambda,\omega^+}$ be the following operator,
$$
H_{\lambda,\omega^+}=-\Delta+W_{\text{per}}+\lambda\sum_{\gamma\in\mathbb{Z}^d}\omega^+f(\cdot-\gamma).
$$
As $H_{\lambda,\omega^+}$ is a $\mathbb{Z}^d$-periodic operator,
the analysis given in subsection \ref{flo} for
$H_{\lambda,\omega^+}$ is still true in the present case. For
$(E_{j}(\lambda,\theta))_{j\in\mathbb{N}^*}$, the Floquet
eigenvalue of $H_{\lambda,\omega^+}$ let us set
 $E_{\lambda,\omega^+}=\inf_{\theta\in \mathbb{T}^*}E_{1}(\lambda,\theta)$.
\begin{thm}\label{ram1}
Assume that the matrix $A(0)$ is negative-definite. Then there
exists $\lambda_{0}>0$ such that for any $\lambda\in
[0,\lambda_0]$ we have:
$$\Pi_{\lambda,0}\Big(H_{\lambda,\omega}-E_{\lambda,\omega^+}\Big)\Pi_{\lambda,0}\
\ \ \text{is\ a \ positive\ operator.}
$$
\end{thm}
The result of Theorem \ref{ram1} can be proved in the same way as
we did for Theorem \ref{hida} in the previous subsection. Indeed,
$H_{\lambda,\omega}$ can be seen as a perturbation of
$H_{\lambda,\omega^+}$ as follow,
$$
H_{\lambda,\omega}=H_{\lambda,\omega^+}+V_{\overline{\omega}}.
$$
With $\displaystyle V_{\overline{\omega}}(\cdot)
=\sum_{\gamma\in\mathbb{Z}^d}\overline{\omega}_{\gamma}f(\cdot-\gamma)$
and for any $\gamma\in\mathbb{Z}^d$
$\overline{\omega}_{\gamma}=\omega_{\gamma}-\omega^+$. Notice that
in this case
$(\overline{\omega}_{\gamma})_{\gamma\in\mathbb{Z}^d}$ is a
family of bounded and negative random variables. \\
using the analogous unitary transformation, one gets that
$H_{\lambda,\omega}^0$ is unitarly equivalent to
$$
h_{\lambda,\omega}^0=h_{\lambda,\omega^+}^0+V_{\lambda,\overline{\omega}}^0.
$$
The lower bound of $h_{\lambda,\omega^+}$ can be derived easily.
As all arguments used to lower bound $V_{\widetilde{\omega}}^0$
remain valid; we lower bound  $V_{\lambda,\overline{\omega}}^0$
using the same computation done in subsection \ref{tak}. \\
So we get that there exist $K_1,K_2>0$ such that
\begin{multline}
\Big|\langle V_{\lambda,\overline{\omega}}^0u,u\rangle
-\sum_{1\leq k,k'\leq m}\sum_{\gamma\in
\mathbb{Z}^d}\overline{\omega}_{\gamma}
\Big(\int_{C_0}a_{\varphi_{\lambda,1},k,k'}(x)dx\Big)(\hat{u_k})(\gamma)
\overline{(\hat{u_{k'}})(\gamma)}\Big|\leq \\
K_1\beta (\omega^+-\omega^-)\sum_{1\leq k\leq m}\sum_{\gamma\in
\mathbb{Z}^d}|(\hat{u_k})(\gamma)|^2+K_2(1+1/\beta)\sum_{1\leq
k\leq m}\langle
\zeta_{k,\lambda}u_k,u_k\rangle_{L^2(\mathbb{T}^*)}.\label{joum1}
\end{multline}
When $A(0)$ is negative-definite, there exists $C<0$ and
$\lambda_0>0$ such that for any $\lambda\in[0,\lambda_0]$, we have
$$A(\lambda)\leq C I_m.$$
As the random variables
$(\overline{\omega}_{\gamma})_{\gamma\in\mathbb{Z}^d}$ are
negative, we get that
\begin{multline}C\sum_{1\leq k\leq
m}\sum_{\gamma\in \mathbb{Z}^d}\overline{\omega}_{\gamma}|(\hat{u_k})(\gamma)|^2\\
\leq \sum_{1\leq k,k'\leq m}\sum_{\gamma\in
\mathbb{Z}^d}\overline{\omega}_{\gamma}\Big(\int_{C_0}
a_{\varphi_{\lambda,1},{k},{k'}}(x)dx\Big)(\hat{u_k})(\gamma)\overline{(\hat{u_{k'}})(\gamma)}.
\label{yas11}\end{multline} This and equation (\ref{joum1}) give
the sought result on the lower bound of $V_{\overline{\omega}}$.
This ends the proof of Theorem \ref{ram1}. \hfill$\Box$\newline
\noindent {\bf{The proof of Lemma \ref{sousou1}:}}\newline
\noindent As $V_{\widetilde{\omega}} $ is
$H_{\lambda,\omega^-}$-relatively bound uniformly on
$\widetilde{\omega}_{\gamma}$, there exists $c>0$ such that for
any $u\in L^{2}(\mathbb{T}^*)$ we have
\begin{eqnarray*}
\| V_{\widetilde{\omega}}{\mathcal{T}}_{\varphi}(u)\| &\leq&
c\Big(\|H_{\lambda,\omega^-}{\mathcal{T}}_{\varphi}(u)\|^2+\|{\mathcal{T}}_{\varphi}(u)\|^2\Big)\\
&\leq
&c\int_{\mathbb{T}^*}\Big(\|H_{\lambda,\omega^-}(\theta)\varphi(\cdot,\theta)\|_{L^{2}(C_0)}^2+\|\varphi(\cdot,\theta)\|^2_{L^2(C_0)}\Big)|u(\theta)|^2d\theta.\\
&\leq c&\|\varphi\|_{H_{\lambda,\omega^-},\infty}^2\cdot
\|u\|^{2}_{L^2(\mathbb{T}^*)}.
\end{eqnarray*}
One computes\\
$\langle
V_{\widetilde{\omega}}^0{\mathcal{T}}_{\varphi_k}(u),{\mathcal{T}}_{\varphi_{k'}}(v)\rangle_{L^{2}(\mathbb{R}^d)}$
\begin{equation*}
\int_{\mathbb{R}^d}
V_{\widetilde{\omega}}(x)\varphi_k(x,\theta_k)\overline{\varphi_{k'}(x,\theta_{k'})}\cdot
\Big(\int_{\mathbb{T}^*}e^{(\theta-\theta_k)x}v\theta)d\theta\cdot
{\overline{\int_{\mathbb{T}^*}e^{(\theta-\theta_{k'})x}u(\theta)d\theta}}\Big)dx
\end{equation*}
\begin{multline*}
\sum_{\gamma\in\mathbb{Z}^d}\widetilde{\omega}_{\gamma}\int_{C_0}
f(x)\varphi_k(x,\theta_k)\overline{\varphi_{k'}(x,\theta_{k'})}\cdot
\Big(\int_{\mathbb{T}^*}e^{(\theta-\theta_k)x}e^{i\gamma\cdot\theta}u(\theta)d\theta\cdot\\
{\overline{\int_{\mathbb{T}^*}e^{(\theta-\theta_{k'})x}e^{i\gamma\cdot
\theta}u(\theta)d\theta}}\Big)dx.
\end{multline*}
Let $$\hat{u}(\gamma)=\int_{\mathbb{T}^*}e^{i\gamma\cdot \theta
}u(\theta)d\theta.$$
 For any
$(x,\theta)\in\mathbb{R}^d\times \mathbb{T}^*$ and $1\leq k\leq
m$, we set
\begin{equation}
g_{k}(x,\theta)=\frac{e^{i(\theta-\theta_k)\cdot
x}-1}{\sqrt{\zeta_{k,\lambda}(\theta)}}.
\end{equation}
As $\theta_{k}(\lambda)$ is the only zero of $\zeta_{k,\lambda}$
and as it is nondegenerate, there exist $C>0$ such that, for
$(x,\theta)\in\mathbb{R}\times \mathbb{R}^d$, and $1\leq k\leq m$,
we have
\begin{equation}
|g_{k}(x,\theta)|\leq C(1+|x|).\label{eq:k}
\end{equation}
We have $\displaystyle
e^{i(\theta-\theta_k)x}=\sqrt{\zeta_{k,\lambda}}g_{k}(x,\theta)+1$.
So using this and expanding $\langle
V_{\widetilde{\omega}}{\mathcal{T}}_{\varphi_k}(u),{\mathcal{T}}_{\varphi_{k'}}(v)\rangle_{L^{2}(\mathbb{R}^d)}$,
we get
\begin{multline}
\langle
V_{\widetilde{\omega}}{\mathcal{T}}_{\varphi_k}(u),{\mathcal{T}}_{\varphi_{k'}}(v)\rangle_{L^{2}(\mathbb{R}^d)}
-\sum_{\gamma\in
\mathbb{Z}^d}\widetilde{\omega}_{\gamma}\Big(\int_{C_0}a_{\varphi,k,k'}(x)dx\Big)\cdot
|u(\gamma)|^2\\
=\sum_{\gamma\in\mathbb{Z}^d}\widetilde{\omega}_{\gamma}\int_{C_0}a_{\varphi,k,k'}(x)\cdot
\Big(\int_{\mathbb{T}^*}e^{i\gamma\cdot
\theta}g_k(x,\theta)\sqrt{\zeta_{k,\lambda}(\theta)}u(\theta)d\theta\cdot\\
\overline{\int_{\mathbb{T}^*} e^{i\gamma\cdot \theta}
g_{k'}(x,\theta)\sqrt{\zeta_{k',\lambda}
(\theta)}v(\theta)d\theta}\Big) dx\\
+\sum_{\gamma\in\mathbb{Z}^d}\widetilde{\omega}_{\gamma}\hat{u}(\gamma)\int_{C_0}a_{\varphi,k,k'}(x)\cdot\overline{\int_{\mathbb{T}^*}e^{i\gamma\cdot
\theta}g_{k'}(x,\theta)\sqrt{\zeta_{k',\lambda}(\theta)}v(\theta)d\theta} dx\\
+\sum_{\gamma\in\mathbb{Z}^d}
\widetilde{\omega}_{\gamma}\overline{\hat{v}(\gamma)}\int_{C_0}a_{\varphi,k,k'}(x)\cdot
\int_{\mathbb{T}^*}e^{i\gamma\cdot
\theta}g_k(x,\theta)\sqrt{\zeta_{k,\lambda}(\theta)}u(\theta)d\theta
dx.\label{er1}
\end{multline}
Now using the fact that the family
$(\widetilde{\omega}_{\gamma})_{\gamma\in\mathbb{Z}^d}$ is bounded, Cauchy-Schwartz inequality and Perseval identity
 and equation (\ref{eq:k}), we get that there exists $C>0$ such that \begin{multline}
\sum_{\gamma\in\mathbb{Z}^d}\widetilde{\omega}_{\gamma}\int_{C_0}a_{\varphi,k,k'}(x)\cdot
\Big(\int_{\mathbb{T}^*}e^{i\gamma\cdot
\theta}g_k(x,\theta)\sqrt{\zeta_{k,\lambda}(\theta)}u(\theta)d\theta\cdot\\
\overline{\int_{\mathbb{T}^*} e^{i\gamma\cdot \theta}
g_{k'}(x,\theta)\sqrt{\zeta_{k',\lambda}
(\theta)}v(\theta)d\theta}\Big) dx\\ \leq
C\langle\zeta_{k,\lambda}u,u\rangle_{L^{2}(\mathbb{T}^*)}+\langle
\zeta_{k'+\lambda}v,v\rangle_{L^{2}(\mathbb{T}^*)}.\label{er2}
\end{multline} And $\beta>0$,
\begin{multline}\sum_{\gamma\in\mathbb{Z}^d}\widetilde{\omega}_{\gamma}\hat{u}(\gamma)\int_{C_0}a_{\varphi,k,k'}(x)\cdot\overline{\int_{\mathbb{T}^*}e^{i\gamma\cdot
\theta}g_{k'}(x,\theta)\sqrt{\zeta_{k',\lambda}(\theta)}v(\theta)d\theta}
dx \leq\\
\beta
\sum_{\gamma\in\mathbb{Z}^d}|\hat{u}(\gamma)|^{2}+\frac{1}{4\beta}\langle
\zeta_{k',\lambda}v,v\rangle_{L^{2}(\mathbb{T}^*)}.\label{er3}
\end{multline}
The same argument gives
\begin{multline}
\sum_{\gamma\in\mathbb{Z}^d}
\widetilde{\omega}_{\gamma}\overline{\hat{v}(\gamma)}\int_{C_0}a_{\varphi,k,k'}(x)\cdot
\int_{\mathbb{T}^*}e^{i\gamma\cdot
\theta}g_k(x,\theta)\sqrt{\zeta_{k,\lambda}(\theta)}u(\theta)d\theta
dx.\\ \leq \beta
\sum_{\gamma\in\mathbb{Z}^d}|\hat{v}(\gamma)|^{2}+\frac{1}{4\beta}\langle
\zeta_{k,\lambda}u,u \rangle_{L^{2}(\mathbb{T}^*)}.\label{er4}
\end{multline}
So from (\ref{er1}), (\ref{er2}), (\ref{er3}) and (\ref{er4}) we
get (\ref{eq03}).\newline The proof of (\ref{eq04}); follows by
changing $\mathcal{T}_{\varphi_{\lambda,1}}(u)$ using
(\ref{wakt1}) and following the same steps as (\ref{eq03}).
\newline This ends the proof of Lemme \ref{sousou1}
\hfill$\Box$ \newline \textit{$\mathbf{Acknowledgements.}$ The
author would like to thank Frederic Klopp for his precious
comments and remarks.}
\bibliographystyle{amsplain}

\end{document}